\documentclass{amsart}
\usepackage{geometry,amssymb,color,hyperref}
\usepackage{graphicx} 
\usepackage{mathtools}

\theoremstyle{plain}

\theoremstyle{definition}

\theoremstyle{remark}
\newtheorem{remark}{Remark}

\newcommand{\IZ}{\mathbb Z}

\title{There exist Steiner systems $S(2,8,225)$ and $S(2,9,289)$}
\author{Ivan Hetman}
\date{September 2025}
\dedicatory{To my grandfather, Stepan Yurchyshyn, who recently turned 92.}

\begin{document}
    \begin{abstract} In this note six Steiner systems $S(2,8,225)$ and four Steiner systems $S(2,9,289)$ are presented. This resolves two of $129$ undecided cases for block designs with block length $8$ and $9$, mentioned in Handbook of Combinatorial Designs.
    \end{abstract}

    \maketitle

    Given three integers $t$, $k$, $v$ such that $2 \le t < k < v$, a Steiner system $S(t, k, v)$ is a set $V$ of cardinality $v$ together with a family $\mathcal B$ of $k$-subsets of $V$ (called \textit{blocks}) such that every $t$-element subset of $V$ is contained in exactly one block. Handbook of Combinatorial Designs \cite[II.3.8-9,11-12]{HoCD} lists $38$ and $91$ numbers $v$ for which the existence of Steiner systems $S(2,8,v)$ and $S(2,9,v)$, respectively, is unknown. In this note we present six Steiner systems $S(2,8,225)$ and four Steiner systems $S(2,9,289)$, thus resolving two of the above $129$ undecided cases.

    We will encode new Steiner systems as difference families preceded by their fingerprints. For definition, properties and examples of difference families in groups, see \cite[VI.16]{HoCD}. The fingerprint of a Steiner system $S(2,k,v)$ is its combinatorial characteristic of geometrical origin, counted as follows. For any distinct elements $a,b\in V$, we denote by $\overline{ab}\in \mathcal B$ a unique block that contains the elements $a,b$. Let $Q$ be the set of all quadruples $(u,x,o,y)$ of distinct elements of $V$ such that the block $\overline{ox}$ contains $u$ but not $y$. For any quadruple $(u,x,o,y)\in Q$, let $N(u,x,o,y)$ be the number of elements $v\in\overline{oy}$ such that the blocks $\overline{uv}$ and $\overline{xy}$ are disjoint.

    \begin{picture}(100,90)(-180,-15)
        \put(0,0){\line(1,0){60}}
        \put(0,0){\line(1,2){30}}
        \put(60,0){\line(-1,2){30}}
        \put(15,30){\line(1,0){30}}

        \put(0,0){\circle*{2}}
        \put(-4,-8){$x$}
        \put(60,0){\circle*{2}}
        \put(59,-8){$y$}
        \put(30,60){\circle*{2}}
        \put(28,63){$o$}
        \put(15,30){\circle*{2}}
        \put(7,28){$u$}
        \put(45,30){\circle*{2}}
        \put(47,28){$v$}
    \end{picture}

    The {\em fingerprint} of the Steiner system $(V,\mathcal B)$ is the function $f$ assigning to each number $n_i$ in  the set $\{n_1,\dots,n_\ell\}\coloneqq\{N(u,x,o,y):(u,x,o,y)\in Q\}$ the cardinality of the set $\{(u,x,o,y)\in Q:N(u,x,o,y)=n_i\}$. We represent the fingerprint $f$ in the form
    \[\{n_1=f(n_1),n_2=f(n_2),\dots,n_\ell=f(n_\ell)\}.\]

    \begin{remark} The definition of the fingerprint ensures that Steiner systems with different fingerprints are not isomorphic. The converse is not true in general, but in most cases, non-isomorphic Steiner systems $S(2,k,v)$ for $k\ge 6$ and $v>k^2$  have different fingerprints. Since the calculation of the fingerprint of a Steiner system $S(2,k,v)$ has computational complexity $O(k^2v^3)$, it is an efficient practical tool for distinguishing non-isomorphic Steiner systems.
    \end{remark}

    Two Steiner systems $S(2,8,225)$ are found as difference families in the Abelian group\break $\mathbb Z_3 \times \mathbb Z_3 \times \mathbb Z_5 \times \mathbb Z_5$ whose elements are labeled by quadruples $abcd$ of numbers
    $a,b\in\{0,1,2\}=\IZ_3$ and $c,d\in \{0,1,2,3,4\}=\IZ_5$:
    \begin{enumerate}
        \item $\{1=16200, 2=263250, 3=2683800, 4=14717700, 5=29170800, 6=18769050\} \newline
        \big\{\{0000, 0001, 0103, 1003, 1210, 1241, 2112, 2144\}$, \newline
        ${}\;\;\{0000, 0002, 0121, 0131, 0222, 0230, 2101, 2201\}$, \newline
        ${}\;\;\{0000, 0011, 1001, 1010, 1233, 2023, 2043, 2233\}$, \newline
        ${}\;\;\{0000, 0012, 1031, 1120, 1142, 2131, 2223, 2244\}\big\}$;
        \item $\{2=202500, 3=2786400, 4=14754150, 5=29211300, 6=18666450\} \newline
        \big\{\{0000, 0001, 0103, 1003, 1210, 1241, 2112, 2144\}$, \newline
        ${}\;\;\{0000, 0002, 0121, 0131, 0222, 0230, 2101, 2201\}$, \newline
        ${}\;\;\{0000, 0011, 1001, 1010, 1233, 2023, 2043, 2233\}$, \newline
        ${}\;\;\{0000, 0012, 1123, 1144, 1231, 2031, 2220, 2242\}\big\}$.
    \end{enumerate}

    The other four Steiner systems $S(2,8,225)$ are found as difference families in the Abelian group $\mathbb Z_5 \times \mathbb Z_5 \times \mathbb Z_9$ whose elements are labeled by triples $abc$ of numbers $a,b\in\{0,1,2,3,4\}$ and $c\in\IZ_9=\{0,1,2,3,4,5,6,7,8\}$:
    \begin{enumerate}
        \item $\{1=16200, 2=229500, 3=3088800, 4=14523300, 5=29484000, 6=18279000\} \newline
        \big\{\{000, 001, 012, 042, 103, 117, 403, 447\}, \; \{000, 002, 146, 227, 245, 315, 337, 416\}$, \newline ${}\;\;\{000, 003, 121, 137, 204, 304, 427, 431\},\; \{000, 004, 021, 031, 213, 238, 328, 343\}\big\}$;
        \item $\{1=16200, 2=270000, 3=3045600, 4=15161850, 5=28800900, 6=18326250\} \newline
        \big\{\{000, 001, 012, 042, 103, 117, 403, 447\}, \; \{000, 002, 145, 224, 246, 316, 334, 415\}$, \newline ${}\;\;\{000, 003, 121, 137, 204, 304, 427, 431\}, \;\{000, 004, 021, 031, 213, 238, 328, 343\}\big\}$;
        \item $\{1=21600, 2=270000, 3=3112200, 4=15187500, 5=28684800, 6=18344700\} \newline
        \big\{\{000, 001, 012, 042, 103, 117, 403, 447\}, \; \{000, 002, 145, 224, 246, 316, 334, 415\}$, \newline ${}\;\;\{000, 003, 121, 137, 204, 304, 427, 431\}, \; \{000, 004, 023, 033, 211, 235, 325, 341\}\big\}$;
        \item $\{2=384750, 3=3450600, 4=15078150, 5=28541700, 6=18165600\} \newline
        \big\{\{000, 001, 012, 042, 103, 117, 403, 447\}, \; \{000, 002, 146, 227, 245, 315, 337, 416\}$, \newline ${}\;\;\{000, 003, 122, 135, 208, 308, 425, 432\}, \; \{000, 004, 021, 031, 213, 238, 328, 343\}\big\}$.
    \end{enumerate}

    Four Steiner systems $S(2,9,289)$ are found as difference families in the Abelian group $\mathbb Z_{17} \times \mathbb Z_{17}$ whose elements are labeled by pairs $xy$ of symbols $x,y\in \mathbb Z_{17}=\{0,1,\dots,9,\mbox{A},\mbox{B},\mbox{C},\mbox{D},\mbox{E},\mbox{F},\mbox{G}\}$:

    \begin{enumerate}
        \item $\{2=135252, 3=2271540, 4=14919336, 5=46511082, 6=65739408, 7=33558102\} \newline
        \big\{\{\mbox{00}, \mbox{01}, \mbox{03}, \mbox{13}, \mbox{22}, \mbox{33}, \mbox{4A}, \mbox{6G}, \mbox{AE}\}, \; \{\mbox{00}, \mbox{04}, 09, \mbox{1F}, 59, \mbox{9D}, \mbox{A8}, \mbox{BG}, \mbox{D9}\}$, \newline ${}\;\;\{00, 06, \mbox{2F}, \mbox{3B}, 58, \mbox{A9}, \mbox{BE}, \mbox{D1}, \mbox{EB}\}, \; \{00, 07, 14, \mbox{3F}, \mbox{5A}, \mbox{62}, \mbox{89}, \mbox{A5}, \mbox{CA}\}\big\}$;
        \item $\{2=145656, 3=2531640, 4=15606000, 5=45964872, 6=65344056, 7=33542496\} \newline
        \big\{\{00, 01, 03, 13, 22, 33, \mbox{4A}, \mbox{6G}, \mbox{AE}\}, \; \{00, 04, 09, \mbox{1F}, 59, \mbox{9D}, \mbox{A8}, \mbox{BG}, \mbox{D9}\}$, \newline ${}\;\;\{00, 06, \mbox{2F}, \mbox{3B}, 58, \mbox{A9}, \mbox{BE}, \mbox{D1}, \mbox{EB}\}, \; \{00, 07, \mbox{5E}, 72, \mbox{9F}, \mbox{B5}, \mbox{CE}, \mbox{E9}, \mbox{G3}\}\big\}$;
        \item $\{2=197676, 3=2548980, 4=14974824, 5=47005272, 6=64924428, 7=33483540\} \newline
        \big\{\{00, 01, 04, 15, 65, \mbox{CC}, \mbox{D1}, \mbox{D6}, \mbox{E1}\}, \; \{00, 02, \mbox{2F}, \mbox{7C}, 89, 92, \mbox{A0}, \mbox{B3}, \mbox{DC}\}$, \newline ${}\;\;\{00, 06, 34, \mbox{9D}, \mbox{A3}, \mbox{AC}, \mbox{C5}, \mbox{CF}, \mbox{E9}\}, \; \{00, 12, 31, \mbox{5F}, 88, \mbox{9G}, \mbox{A8}, \mbox{B4}, \mbox{EA}\}\big\}$;
        \item $\{2=239292, 3=2522970, 4=15564384, 5=46443456, 6=65125572, 7=33239046\} \newline
        \big\{\{00, 01, 04, 15, 65, \mbox{CC}, \mbox{D1}, \mbox{D6}, \mbox{E1}\}, \; \{00, 02, \mbox{2F}, \mbox{7C}, 89, 92, \mbox{A0}, \mbox{B3}, \mbox{DC}\}$, \newline ${}\;\;\{00, 06, 34, \mbox{9D}, \mbox{A3}, \mbox{AC}, \mbox{C5}, \mbox{CF}, \mbox{E9}\}, \; \{00, 12, 49, \mbox{7F}, \mbox{8B}, 93, \mbox{AB}, \mbox{D4}, \mbox{F1}\}\big\}$.
    \end{enumerate}

    \section*{Acknowledgement}

    I would like to express my sincere thanks to Taras Banakh and Alex Ravsky for helping me to direct my search and ultimately discover the long-awaited Steiner systems $S(2,8,225)$ and $S(2,9,289)$.

\end{document}